\documentclass[10pt,letterpaper]{amsart} % changed to 10pt
\usepackage{epic,eepic, euler, palatino, xypic, microtype}
 \newlength{\baseunit}               % the basic unit length
 \newcount{\numlines}                % depth of picture (in number of lines)
\newcommand{\tpoint}[1]{\vspace{3mm}\par \noindent \refstepcounter{subsection}{\bf \thesubsection.} 
  {\em #1. ---} }

\newcommand{\proj}{\mathbb P}
\newcommand{\oh}{{\mathcal{O}}}

\newcommand{\si}{\sigma}

\newcommand{\Sec}{\operatorname{Sec}}
\newcommand{\Sing}{\operatorname{Sing}}

\newcommand{\Sym}{\operatorname{Sym}}
\newcommand{\Aut}{\operatorname{Aut}}
\newcommand{\Z}{\mathbb{Z}}

\newcommand{\vw}{\vec{w}}

\begin{document}
\pagestyle{plain}

\title{\large{The relations among invariants of points on the projective line}}
\author{Ben Howard, John Millson, Andrew Snowden, and Ravi Vakil}

\thanks{*B.\ Howard was supported by NSF fellowship DMS-0703674.
J.\ Millson was supported by NSF grant DMS-0405606, NSF FRG grant DMS-0554254 and the Simons Foundation.  R.\ Vakil was supported by NSF grant DMS-0801196}

\date{Friday, June 12, 2009.}
\begin{abstract}
 We consider the ring of invariants of $n$ points on the projective
  line.  The space $(\proj^1)^n / \! / \rm{PGL}_2$ is perhaps the
  first nontrivial example of a Geometry Invariant Theory quotient.
  The construction depends on the weighting of the $n$ points.  Kempe
  \cite{Kempe} discovered a beautiful set of generators 
 (at least in the case of unit weights) in 1894.
  We describe the full ideal of relations for
  all possible weightings.  In some sense, there is only one equation,
  which is quadric except for the classical case of the Segre cubic
  primal, for $n=6$ and weight $1^6$.  The cases of up to $6$ points
  are long known to relate to beautiful familiar geometry.  The case of $8$
  points turns out to be richer still.
\end{abstract}

\maketitle
%\tableofcontents

{\parskip=12pt % closing bracket is just before the bibliography 

\section{Introduction}

We consider the ring of invariants of $n$ points
on the projective line, and the GIT
quotient $(\proj^1)^n / \! / \rm{PGL}_2$.  The quotient depends on a choice
of $n$ weights $\vw := (w_1,\dots, w_n)  \in (\Z^+)^n$. The quotient
is given by 
$$(\proj^1)^n \dashrightarrow (\proj^1)^n  / \! /_{\vw}  \rm{PGL}_2 := \operatorname{Proj} \left( \bigoplus_k R_{k \vec{w}} \right) $$ 
where $R_{\vec{v}} = \Gamma( (\proj^1)^n, \oh( v_1, \dots, v_n))^{\rm{PGL}_2}.$
%Denote the $k$th graded piece of the ring on the right $R^{\vw}_k$.
%If $\vw = 1^n$, we right $R^{(n)}_k$ for convenience.
Small cases ($n \leq 6$)  yield familiar beautiful geometry, and we refer in
particular to \cite{do} for a masterful discussion,
including of the history.  In
these examples, we take all weights to be $1$, and work over a field.
With sufficient care, this entire discussion applies over $\Z$.

The case $n=4$ gives the cross ratio $(\proj^1)^4 \dashrightarrow
\overline{\mathcal{M}}_{0,4} \cong \proj^1$.  From the perspective of
this project, the cross ratio is best understood as a rational map not
to $\proj^1$, but to the line in $\proj^2$ where the three coordinates
sum to $0$; the $S_4$-action is more transparent in this way, as we
will see shortly.

The case $n=5$ yields the quintic del Pezzo surface
$(\proj^1)^5 \dashrightarrow \overline{\mathcal{M}}_{0,5} \hookrightarrow \proj^5$. 
The map is given by the degree $2$ invariants; there are
no degree $1$ invariants.  The quotient is cut out by 5 quadrics.

The case $n=6$ is particularly beautiful,
see Figure~\ref{figsix}. The ring is generated in degree $1$, and this
piece has dimension $5$, so the quotient threefold, the {\em Segre
  cubic}, is naturally a hypersurface in $\proj^4$.  The equations
(for characteristic not $3$) are cleanest written in six variables
$x_1$, \dots, $x_6$:
$$
x_1+x_2+x_3+x_4+x_5+x_6 = 
x_1^3+x_2^3+x_3^3+x_4^3+x_5^3+x_6^3 = 0.
$$
(This description was given by Joubert
in 1867, and in 1911 Coble gave an
invariant-theoretical interpretation of the Joubert identities \cite{coble}, see
also Coble's book \cite[Ch.~ III]{coblebook}.)  There
are two obvious $S_6$ actions on the invariants, the first via
permutations of the points, and the second by permutations of the
variables $x_i$.  One might hope that this would arise from a correspondence between
the points and the variables, but the actions are related by the outer
automorphism of $S_6$.
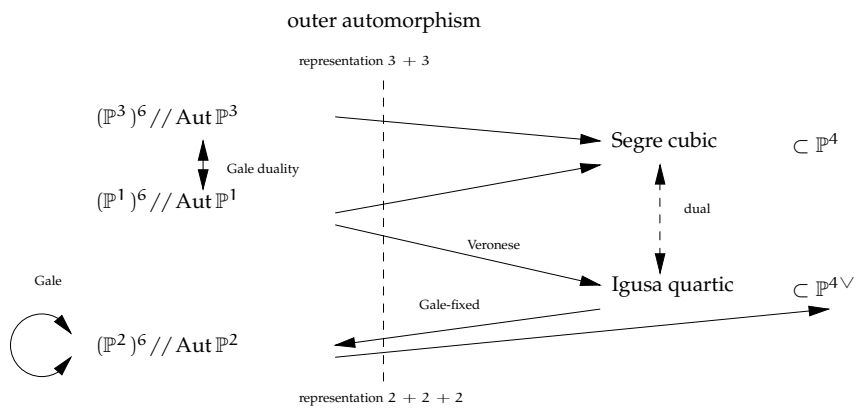
\begin{figure}
\begin{center}
\setlength{\unitlength}{0.00083333in}
\begingroup\makeatletter\ifx\SetFigFont\undefined%
\gdef\SetFigFont#1#2#3#4#5{%
  \reset@font\fontsize{#1}{#2pt}%
  \fontfamily{#3}\fontseries{#4}\fontshape{#5}%
  \selectfont}%
\fi\endgroup%
{\renewcommand{\dashlinestretch}{30}
\begin{picture}(5121,2463)(0,-10)
\put(1809,27){\makebox(0,0)[lb]{{\SetFigFont{5}{6.0}{\rmdefault}{\mddefault}{\updefault}representation $2+2+2$}}}
\path(2034,1827)(3684,1677)
\blacken\path(3561.777,1657.987)(3684.000,1677.000)(3567.209,1717.741)(3561.777,1657.987)
\path(2034,1227)(3684,1527)
\blacken\path(3571.302,1476.018)(3684.000,1527.000)(3560.569,1535.050)(3571.302,1476.018)
\path(2034,1152)(3684,777)
\blacken\path(3560.335,774.341)(3684.000,777.000)(3573.633,832.849)(3560.335,774.341)
\path(2034,327)(5109,627)
\blacken\path(4992.480,585.490)(5109.000,627.000)(4986.654,645.206)(4992.480,585.490)
\blacken\path(1239.000,1557.000)(1209.000,1677.000)(1179.000,1557.000)(1239.000,1557.000)
\path(1209,1677)(1209,1377)
\blacken\path(1179.000,1497.000)(1209.000,1377.000)(1239.000,1497.000)(1179.000,1497.000)
\blacken\path(4089.000,1407.000)(4059.000,1527.000)(4029.000,1407.000)(4089.000,1407.000)
\dashline{60.000}(4059,1527)(4059,852)
\blacken\path(4029.000,972.000)(4059.000,852.000)(4089.000,972.000)(4029.000,972.000)
\dashline{60.000}(2334,2052)(2334,177)
\path(3684,627)(2034,402)
\blacken\path(2148.846,447.938)(2034.000,402.000)(2156.953,388.489)(2148.846,447.938)
\put(534,1752){\makebox(0,0)[lb]{{\SetFigFont{8}{9.6}{\rmdefault}{\mddefault}{\updefault}$(\proj^3)^6 // \Aut \proj^3$}}}
\put(534,1227){\makebox(0,0)[lb]{{\SetFigFont{8}{9.6}{\rmdefault}{\mddefault}{\updefault}$(\proj^1)^6// \Aut \proj^1$}}}
\put(3759,1602){\makebox(0,0)[lb]{{\SetFigFont{8}{9.6}{\rmdefault}{\mddefault}{\updefault}Segre cubic}}}
\put(1359,1452){\makebox(0,0)[lb]{{\SetFigFont{5}{6.0}{\rmdefault}{\mddefault}{\updefault}Gale duality}}}
\put(3759,702){\makebox(0,0)[lb]{{\SetFigFont{8}{9.6}{\rmdefault}{\mddefault}{\updefault}Igusa quartic}}}
\put(4884,702){\makebox(0,0)[lb]{{\SetFigFont{8}{9.6}{\rmdefault}{\mddefault}{\updefault}$\subset \proj^{4 \vee}$}}}
\put(4884,1602){\makebox(0,0)[lb]{{\SetFigFont{8}{9.6}{\rmdefault}{\mddefault}{\updefault}$\subset \proj^4$}}}
\put(534,327){\makebox(0,0)[lb]{{\SetFigFont{8}{9.6}{\rmdefault}{\mddefault}{\updefault}$(\proj^2)^6 // \Aut \proj^2$}}}
\put(2859,1002){\makebox(0,0)[lb]{{\SetFigFont{5}{6.0}{\rmdefault}{\mddefault}{\updefault}Veronese}}}
\put(159,777){\makebox(0,0)[lb]{{\SetFigFont{5}{6.0}{\rmdefault}{\mddefault}{\updefault}Gale}}}
\put(4209,1227){\makebox(0,0)[lb]{{\SetFigFont{5}{6.0}{\rmdefault}{\mddefault}{\updefault}dual}}}
\put(2559,627){\makebox(0,0)[lb]{{\SetFigFont{5}{6.0}{\rmdefault}{\mddefault}{\updefault}Gale-fixed}}}
\put(1809,2127){\makebox(0,0)[lb]{{\SetFigFont{5}{6.0}{\rmdefault}{\mddefault}{\updefault}representation $3+3$}}}
\put(1734,2352){\makebox(0,0)[lb]{{\SetFigFont{8}{9.6}{\rmdefault}{\mddefault}{\updefault}outer automorphism}}}
\put(204.000,402.000){\arc{390.000}{0.3948}{5.8884}}
\blacken\path(282.729,548.022)(384.000,477.000)(328.066,587.324)(282.729,548.022)
\blacken\path(328.066,216.676)(384.000,327.000)(282.729,255.978)(328.066,216.676)
\end{picture}
}
\end{center}
\caption{The classical geometry of six points in projective space}
\label{figsix}
\end{figure}

The quotient $(\proj^3)^6/ \! / \Aut( \proj^3)$ is canonically
isomorphic to $(\proj^1)^6/ \! / \Aut( \proj^1)$, via the Gale
transform.  The quotient $(\proj^2)^6 / \! / \Aut (\proj^2)$ is a
double cover of $\proj^4$, branched over the Igusa quartic
hypersurface, given by the equations
$$
w_1+ \cdots + w_6 = 4 (w_1^4 + \cdots + w_6^4) -  
(w_1^2 + \cdots + w_6^2)^2 = 0.
$$
The points and the variables are again related by the outer automorphism of
$S_6$.   Gale duality exchanges the two sheets of the double cover, and thus
the Gale-fixed points,  where the six points lie on a conic, 
correspond to the Igusa quartic.  The Igusa quartic is dual
(in the sense of classical projective geometry)  to the Segre cubic, and
the map on moduli sends six points on $\proj^1$ to six points on the conic.

Our motivating question is the following: {\em is there similarly rich
 structure in the case of more points}?  In \S \ref{sectioneight}, we
describe a structure for eight points parallel to, and in some sense
generalizing, that of six points.  Further discussion (and proofs)
will be given in \cite{hmsv8}.  Although the results are geometric,
the proofs are essentially representation theory.  In \S \ref{big}, we
describe the generators for the ideal of relations for $n$ points in
general.  They are  ``inherited'' from the $n=8$ case.
This completes the program initiated in \cite{hmsv1}, and
details and proofs will be given in \cite{hmsv2}.

\section{Eight points}

\label{sectioneight}
The quotient $M_8 := (\proj^1)^8 / \! / \Aut(\proj^1)$ naturally lies in
$\proj^{13}$.  (The graded ring is generated in degree $1$, and its
degree $1$ piece has rank $14$.)  A number of authors have shown by
computer calculation that the ideal of relations in characteristic $0$
is generated by $14$ quadrics (Koike, Kondo, Freitag, Salvati Manni,
Maclagan, ...).  We describe the structure of the quadrics more
directly, motivated by a suggestion of Dolgachev.  See
Figure~\ref{figureeight}.

The symmetric group acts on the graded ring.  There is a unique
skew-invariant cubic, and it lies in the ideal of $M_8$ because it is
divisible by the Vandermonde polynomial upon pullback to
$(\proj^1)^8$.  (Sam Grushevsky and Riccardo Salvati Manni have pointed out to us that this fact
that it contains $M_8$ readily follows in the language of theta functions,
see \cite[\S 5]{fsm}.)  In fact more is true: $M_8$ is the
singular locus of the skew cubic, and furthermore the cone over $M_8$ is
scheme-theoretically the singular locus of affine cone of the cubic.
Thus the $14$ quadrics are the $14$ partial derivatives of the skew cubic (and
in characteristic $3$, the cubic itself is a necessary generator of
the ideal).  We emphasize that no computer verification is required.

The dual to the skew cubic has surprisingly low degree --- it is a skew quintic
in $\proj^{13}$, whose singular locus $N'_8$ has dimension $9$.  The
moduli space $N_8$ of $8$ points in $\proj^3$ is a double cover of
this singular locus.  The sheets of the double cover are exchanged by
Gale-duality, which sends $8$ points in $\proj^3$ to $8$ points in
$\proj^3$ (both up to projective equivalence).

By Bezout's theorem, the secant variety to $M_8$ is contained in the
skew cubic; it is a divisor.  The duality birational map $f$ from the cubic
to the quintic blows down this divisor: given a secant line to $M_8$
meeting $M_8$ at points $p$ and $q$, the $14$ quadrics (which give the
duality map) vanish at $p$ and $q$, and thus are scalar multiples of
each other.  Hence we have identified one of the two dimensions of the
$\Sec M_8$ contracted by $f$.  We  now describe the other dimension
contracted, by lifting $f: \Sec(M_8) \dashrightarrow N'_8$ to
$\Sec(M_8) \dashrightarrow N_8$. Fix two distinct points of $M_8$, and
hence a secant line.  From this data of a pair of octuples of points
on $\proj^1$, we obtain an octuple of points on $\proj^1 \times \proj^1$,
which via the Segre embedding yields $8$ points in $\proj^3$ (up to
projective equivalence).  Conversely, given $8$ general points in
$\proj^3$, we can find a smooth quadric through them (indeed, a
one-parameter of smooth quadrics), yielding an unordered pair of
octuples of points on $\proj^1$.  (In fact this construction is used to show
that the skew quintic is dual to the skew cubic.)

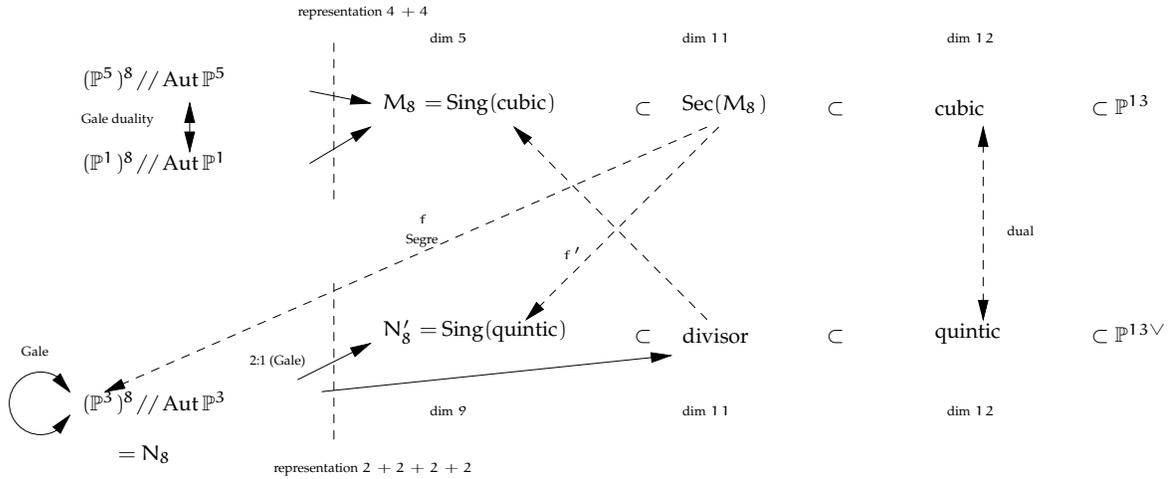
\begin{figure}
\begin{center}
\setlength{\unitlength}{0.00083333in}
\begingroup\makeatletter\ifx\SetFigFont\undefined%
\gdef\SetFigFont#1#2#3#4#5{%
  \reset@font\fontsize{#1}{#2pt}%
  \fontfamily{#3}\fontseries{#4}\fontshape{#5}%
  \selectfont}%
\fi\endgroup%
{\renewcommand{\dashlinestretch}{30}
\begin{picture}(6759,2964)(0,-10)
\put(4209,402){\makebox(0,0)[lb]{{\SetFigFont{5}{6.0}{\rmdefault}{\mddefault}{\updefault}$\dim 11$}}}
\path(1884,2427)(2259,2352)
\blacken\path(2135.447,2346.117)(2259.000,2352.000)(2147.214,2404.951)(2135.447,2346.117)
\path(1884,1977)(2259,2202)
\blacken\path(2171.536,2114.536)(2259.000,2202.000)(2140.666,2165.985)(2171.536,2114.536)
\blacken\path(2165.085,771.502)(2259.000,852.000)(2138.252,825.167)(2165.085,771.502)
\path(2259,852)(1809,627)
\path(1959,552)(4134,777)
\blacken\path(4017.724,734.811)(4134.000,777.000)(4011.550,794.493)(4017.724,734.811)
\blacken\path(6114.000,2082.000)(6084.000,2202.000)(6054.000,2082.000)(6114.000,2082.000)
\dashline{60.000}(6084,2202)(6084,1002)
\blacken\path(6054.000,1122.000)(6084.000,1002.000)(6114.000,1122.000)(6054.000,1122.000)
\dashline{60.000}(4359,1002)(3159,2202)
\blacken\path(3265.066,2138.360)(3159.000,2202.000)(3222.640,2095.934)(3265.066,2138.360)
\dashline{60.000}(4434,2202)(3234,1002)
\blacken\path(3297.640,1108.066)(3234.000,1002.000)(3340.066,1065.640)(3297.640,1108.066)
\dashline{60.000}(4359,2202)(609,552)
\blacken\path(706.756,627.788)(609.000,552.000)(730.920,572.869)(706.756,627.788)
\dashline{60.000}(2034,2727)(2034,1752)
\blacken\path(1164.000,2232.000)(1134.000,2352.000)(1104.000,2232.000)(1164.000,2232.000)
\path(1134,2352)(1134,2052)
\blacken\path(1104.000,2172.000)(1134.000,2052.000)(1164.000,2172.000)(1104.000,2172.000)
\dashline{60.000}(2034,1227)(2034,252)
\put(2334,2277){\makebox(0,0)[lb]{{\SetFigFont{8}{9.6}{\rmdefault}{\mddefault}{\updefault}$M_8 = \Sing(\text{cubic})$ }}}
\put(2334,852){\makebox(0,0)[lb]{{\SetFigFont{8}{9.6}{\rmdefault}{\mddefault}{\updefault}$N'_8 = \Sing(\text{quintic})$}}}
\put(4209,2277){\makebox(0,0)[lb]{{\SetFigFont{8}{9.6}{\rmdefault}{\mddefault}{\updefault}$\Sec(M_8)$}}}
\put(3909,2277){\makebox(0,0)[lb]{{\SetFigFont{8}{9.6}{\rmdefault}{\mddefault}{\updefault}$\subset$}}}
\put(5784,2277){\makebox(0,0)[lb]{{\SetFigFont{8}{9.6}{\rmdefault}{\mddefault}{\updefault}$\text{cubic}$}}}
\put(3909,852){\makebox(0,0)[lb]{{\SetFigFont{8}{9.6}{\rmdefault}{\mddefault}{\updefault}$\subset$}}}
\put(5109,852){\makebox(0,0)[lb]{{\SetFigFont{8}{9.6}{\rmdefault}{\mddefault}{\updefault}$\subset$}}}
\put(5784,852){\makebox(0,0)[lb]{{\SetFigFont{8}{9.6}{\rmdefault}{\mddefault}{\updefault}$\text{quintic}$}}}
\put(6234,1527){\makebox(0,0)[lb]{{\SetFigFont{5}{6.0}{\rmdefault}{\mddefault}{\updefault}dual}}}
\put(2634,2727){\makebox(0,0)[lb]{{\SetFigFont{5}{6.0}{\rmdefault}{\mddefault}{\updefault}$\dim 5$}}}
\put(4209,2727){\makebox(0,0)[lb]{{\SetFigFont{5}{6.0}{\rmdefault}{\mddefault}{\updefault}$\dim 11$}}}
\put(5859,2727){\makebox(0,0)[lb]{{\SetFigFont{5}{6.0}{\rmdefault}{\mddefault}{\updefault}$\dim 12$}}}
\put(5859,402){\makebox(0,0)[lb]{{\SetFigFont{5}{6.0}{\rmdefault}{\mddefault}{\updefault}$\dim 12$}}}
\put(459,1902){\makebox(0,0)[lb]{{\SetFigFont{8}{9.6}{\rmdefault}{\mddefault}{\updefault}$(\proj^1)^8// \Aut \proj^1$}}}
\put(459,402){\makebox(0,0)[lb]{{\SetFigFont{8}{9.6}{\rmdefault}{\mddefault}{\updefault}$(\proj^3)^8 // \Aut \proj^3$}}}
\put(84,777){\makebox(0,0)[lb]{{\SetFigFont{5}{6.0}{\rmdefault}{\mddefault}{\updefault}Gale}}}
\put(459,2427){\makebox(0,0)[lb]{{\SetFigFont{8}{9.6}{\rmdefault}{\mddefault}{\updefault}$(\proj^5)^8 // \Aut \proj^5$}}}
\put(1509,702){\makebox(0,0)[lb]{{\SetFigFont{5}{6.0}{\rmdefault}{\mddefault}{\updefault}2:1 (Gale)}}}
\put(1809,2877){\makebox(0,0)[lb]{{\SetFigFont{5}{6.0}{\rmdefault}{\mddefault}{\updefault}representation $4+4$}}}
\put(1659,27){\makebox(0,0)[lb]{{\SetFigFont{5}{6.0}{\rmdefault}{\mddefault}{\updefault}representation $2+2+2+2$}}}
\put(6759,2277){\makebox(0,0)[lb]{{\SetFigFont{8}{9.6}{\rmdefault}{\mddefault}{\updefault}$\subset \proj^{13}$}}}
\put(6759,852){\makebox(0,0)[lb]{{\SetFigFont{8}{9.6}{\rmdefault}{\mddefault}{\updefault}$\subset \proj^{13 \vee}$}}}
\put(2484,1452){\makebox(0,0)[lb]{{\SetFigFont{5}{6.0}{\rmdefault}{\mddefault}{\updefault}Segre}}}
\put(5109,2277){\makebox(0,0)[lb]{{\SetFigFont{8}{9.6}{\rmdefault}{\mddefault}{\updefault}$\subset$}}}
\put(4209,852){\makebox(0,0)[lb]{{\SetFigFont{8}{9.6}{\rmdefault}{\mddefault}{\updefault}divisor}}}
\put(684,102){\makebox(0,0)[lb]{{\SetFigFont{8}{9.6}{\rmdefault}{\mddefault}{\updefault}$= N_8$}}}
\put(2559,1602){\makebox(0,0)[lb]{{\SetFigFont{5}{6.0}{\rmdefault}{\mddefault}{\updefault}$f$}}}
\put(3459,1377){\makebox(0,0)[lb]{{\SetFigFont{5}{6.0}{\rmdefault}{\mddefault}{\updefault}$f'$}}}
\put(459,2202){\makebox(0,0)[lb]{{\SetFigFont{5}{6.0}{\rmdefault}{\mddefault}{\updefault}Gale duality}}}
\put(2634,402){\makebox(0,0)[lb]{{\SetFigFont{5}{6.0}{\rmdefault}{\mddefault}{\updefault}$\dim 9$}}}
\put(204.000,477.000){\arc{390.000}{0.3948}{5.8884}}
\blacken\path(282.729,623.022)(384.000,552.000)(328.066,662.324)(282.729,623.022)
\blacken\path(328.066,291.676)(384.000,402.000)(282.729,330.978)(328.066,291.676)
\end{picture}
}
\end{center}
\caption{Relations among moduli spaces of eight points in projective space}
\label{figureeight}
\end{figure}

Even more structure may be present. The trisecant variety to $N'_8$ is contained
in the skew quintic, and the quadrisecant variety to $N'_8$ is contained in the divisor.
A dimension estimate shows that this is quite unusual; both should easily fill out the
full $\proj^{13 \vee}$.  One might expect that both containments are equality, although
we have not yet shown this.

This generalizes the $n=6$ story in a number of ways.  For example,
the analogue of the cubic for $n=8$ is the Segre cubic for $n=6$.
(The analogue for $n=4$ is the union of the three boundary points.)  Also, the similarity between
Figures~\ref{figsix} and~\ref{figureeight} is not coincidental: one
can describe a fibration over the constructions of Figure \ref{figsix}
in the boundary of Figure~\ref{figureeight}, commuting with all
dualities.  For example, if we consider octuples of points in
$\proj^3$ where the two points coincide, projecting from those
two points yields six points in $\proj^2$, and the Gale dualities in the two 
figures correspond.

\section{The general case}
\label{big}

We describe the invariants in terms of a {\em graphical algebra}.
To a directed graph $\Gamma$ (with no loops) on $n$ ordered vertices (in bijection
with the $n$ points), we associate
$$
\prod_{\vec{ab} \in \Gamma} (x_a y_b - y_a x_b),$$
an invariant element of $\oh(\vec{v})$, where $\vec{v}$
is the $n$-tuple of valences of the vertices.
The degree $\vw$ invariants are generated (as a vector space or module) by these elements.

This description can be used to show that the ring of invariants
for any $\vw$ is generated in degree $1$.  In the unit weight case,
this is Kempe's Theorem \cite{Kempe}.    

{\em Remark.}  Weyl's
theorems on rings of invariants for a group representation in \cite{weyl} are of the
form: {\em First Main Theorem:} describe generators of the ring of
invariants.  {\em Second Main Theorem:} describe relations among the
generators.  One of his main results is for the symplectic group
acting diagonally on the direct sum of $n$ copies of the standard
representation (\cite[Thm.~ 6.1.A and 6.1.B]{weyl} are the first and second main theorems
respectively); for the case of $\rm{SL}_2$, we obtain the
affine cone over the Grassmannian $Gr(2,n)$.  The first main theorem gave generators for
the invariants, the Pl\"ucker coordinates.  The second main theorem
gave the relations, the Pl\"ucker relations.  Kempe \cite{Kempe} 
proved the first main theorem when the
direct sum is replaced by the tensor product.  Theorem~\ref{mainthmhere}
(or more correctly the main theorem of \cite{hmsv2}) is the second main theorem.

We make a series of observations about this graphical algebra.

\begin{figure}
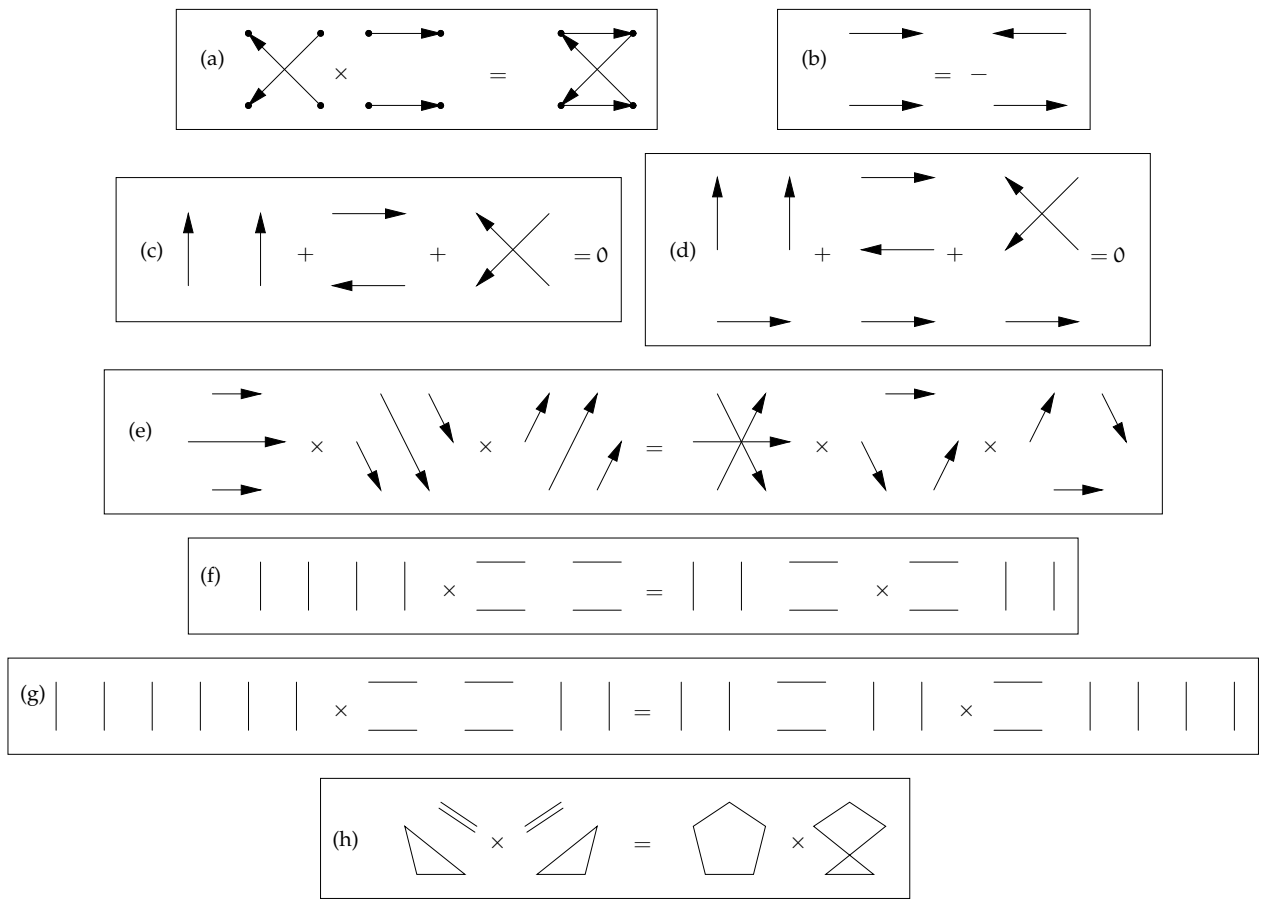

\begin{center}
\include{graphical}
\end{center}
\caption{Relations in the graphical algebra}
\label{fun}
\end{figure}

{\em Multiplication.}  Multiplication of (elements associated to)
graphs is by superposition.  (See for example Figure~\ref{fun}(a).
The vertex labels $1$ through $4$ are omitted for simplicity.  In
later figures, even the vertices will be left implicit.)

{\em Sign (linear) relations.}  Changing the orientation of a single edge
changes the sign of the invariant (e.g.\ Figure~\ref{fun}(b)).

{\em Pl\"ucker (linear) relation.}  Direct calculation shows
the relation of Figure~\ref{fun}(c).

{\em Bigger relations from smaller ones.} The ``four-point'' Pl\"ucker
relation immediately ``extends'' to relations among more points, e.g.\
Figure~\ref{fun}(d) for 6 points.  Any relation may be extended in
this way.  For example, the sign relation in general should be seen as
an extension of the two-point sign relation.

{\em Remark.}  The sign and (extended) Pl\"ucker relations generate
all the linear relations.  They can be used to show that the
``non-crossing graphs'' (with no pairs of edges crossing) with all
edges oriented ``upwards'' (the arrow points toward the higher-numbered label) form a basis --- the graphical form of the
straightening algorithm.  But breaking symmetry obscures the
structure of the ring.

{\em The Segre cubic.}  The relation  of Figure~\ref{fun}(e) is patently true:
the superposition of the three graphs on the left is the same as that of the three
graphs on the right.  This is a cubic relation on the six point space.
It turns out to be nonzero, and is thus necessarily the Segre cubic relation.  
Of course, all that matters about the orientations of the edges is that they are the same
on the both sides of the equation.

{\em The skew cubic on eight points.}  One can describe the Segre
cubic differently: if $\Gamma$ is any $1$-regular graph, consider
$\sum_{\si \in S_6} \operatorname{sgn}(\si) \si(\Gamma^3)$.   Replacing
$6$ by $8$ yields the skew cubic of \S
\ref{sectioneight}.  Replacing $6$ by (even) $n>8$ yields $0$, but there are still nontrivial analogues for arbitrary $n$.

{\em A simple (binomial) quadric on eight points.}
Figure~\ref{fun}(f) gives an obvious relation on $8$ points.  The arrowheads are omitted for
simplicity; they should be chosen consistently on both sides, as in Figure~\ref{fun}(e).

{\em Simple quadrics for at least eight points} are obtained by ``extending'' the eight-point
relations, e.g. Figure~\ref{fun}(g) is the extension to $12$ points, where
the same two edges are added to each graph in Figure~\ref{fun}(f).

We may now state the main theorem of \cite{hmsv2}, in a special case.

\tpoint{Main Theorem for the $n$ even ``unit weight'' case $\vw = 1^n$}
{\em If $n \neq 6$, the simple quadrics (i.e.\ the $S_n$-orbit of the
quadric above) generate the ideal of relations.} \label{mainthmhere}

%This implies that the quotient is explicitly given as  a linear slice
%of a toric variety, which is useful for applications.  Note that it is
%not apparent that the analogue of the skew cubic for $n$ points lies
%in the ideal of the extended simple quadrics, so the theorem cannot be
%simply combinatorial manipulation.

By \cite[Thm. 1.2]{hmsv1}, the arbitrary weight case readily
reduces to the ``unit weight'' case $\vw = 1^n$ ($n$ even), so this
solves the problem for arbitrary weight.  As an example, an explicit
description of the quadrics in the del Pezzo case of five points are as the five
rotations of the patently true relation in Figure~\ref{fun}(h).

%{\em One-sentence explanation:} There is essentially only one equation
%up to symmetry.  There is precisely one equation in the
%``unit weight'' case, and in general, the breaking of symmetry yields
%a number of equations.  This one equation is the Segre cubic if
%$\vw=111111$, and in every other case, the relation is quadric.

{\bf Themes of proof:} 
(i) Use symmetry when necessary.  (ii) Break symmetry when necessary.
(iii) The eight point case is the base of the induction.

More precisely, we use a construction of Speyer-Sturmfels
\cite{SpeyerSturmfels} to degenerate the quotient into a toric
variety, and we show the toric variety is cut out in degree two and
three.  We identify the cubics, and lift them to explicit cubic
relations for the original (quotient) variety.  We show by
representation theory (and elementary combinatorics) that they lie in
the ideal cut out by the quadrics.  Then we use representation theory
to see that the quadrics are generated as a module by the simple
quadrics defined above.

The degree $1$ part $R_1$ of the ring carries the irreducible $S_n$
representation corresponding to $n/2+n/2$.  A key fact is the
following: $\Sym^2 R_1$ is multiplicity free, and consists of irreducibles
corresponding to  all partitions with at most four
parts, all even.  $R_2$ corresponds to those partitions with at most three
parts (all even).  Thus the degree $2$ part of the ideal corresponds to
partitions with precisely four parts, all even.  From this perspective,
 the $n=8$
case is clearly special.  Less obviously, even the form of the simple
binomial quadrics is suggested by the representation theory.

{\bf Conclusion.} We have thus answered our
motivating question: there is sufficient structure in the general
case that we can describe (generators of) the relations completely.
The structure is inherited from the case of $n=8$, where it is a consequence
of exceptional geometry.

{\bf Future prospects.} 
In our situation, we have a family of graded rings where relations in
one case extend to relations in larger cases.  One might hope that
analogues of our results hold in a more general abstract setting.  The
third author has a precise algebraic conjecture implying this in a
wide variety of circumstances generalizing this instance.  In the
particular case of $n$ points in $\proj^1$, this ``coherence
conjecture'' would imply that the relations ``stabilize'' after a
certain number of points, but would not predict an $n$ for which it
would happen.

{\bf Acknowledgments.}  We are indebted to a large number of mathematicians
for helpful comments.  We single out in particular Igor Dolgachev.

} % end of parskip; it started just before the introduction
\end{document}